\newif\ifpictures
\picturestrue

\documentclass[12pt]{amsart}
\usepackage{amssymb}
\usepackage{epsf,epsfig}

\usepackage{enumerate}

\numberwithin{equation}{section}

\newtheorem{thm}{Theorem}
\newtheorem{prop}[thm]{Proposition}
\newtheorem{lemma}[thm]{Lemma}
\newtheorem{cor}[thm]{Corollary}

\newtheorem{example}[thm]{Example}
\newtheorem{remark1}[thm]{Remark}

\newcounter{FNC}[page]
\def\newfootnote#1{{\addtocounter{FNC}{2}$^\fnsymbol{FNC}$%
     \let\thefootnote\relax\footnotetext{$^\fnsymbol{FNC}$#1}}}

\newcommand{\N}{\mathbb{N}}

\newcommand{\R}{\mathbb{R}}
\newcommand{\E}{\mathbb{E}}
\newcommand{\B}{\mathbb{B}}

\DeclareMathOperator{\conv}{conv}
\DeclareMathOperator{\dist}{dist}
\DeclareMathOperator{\lin}{lin}
\DeclareMathOperator{\aff}{aff}

\title[Radii minimal projections
  and optimization of symmetric
  polynomials]{Radii minimal projections of polytopes and constrained 
  optimization of symmetric polynomials}

\author{Ren\'{e} Brandenberg}
\address{Ren\'{e} Brandenberg \\
        Zentrum Mathematik\\
        Technische Universit\"at M\"unchen\\
        Boltzmannstr.~3\\
        D--85747 Gar\-ching bei M\"unchen}
\email{brandenb@ma.tum.de}
\urladdr{http://www-m9.ma.tum.de/\~{}brandenb/}

\author{Thorsten Theobald}
\address{Thorsten Theobald \\
        Institut f\"ur Mathematik \\ MA 6--2 \\
        Technische Universit\"at Berlin \\
        Stra{\ss}e des 17.~Juni 136 \\
        D--10623 Berlin}
\email{theobald@math.tu-berlin.de}
\urladdr{http://www.math.tu-berlin.de/\~{}theobald/}

\date{\today}
\keywords{Polytope, projection, outer radius,
  enclosing cylinder, regular simplex, polynomial optimization, symmetric polynomials.\\
\indent
\emph{2000 Mathematics Subject Classification}
  51N20, 52B12, 52B55, 68U05, 90C26.
}

\begin{document}

\begin{abstract}
  We provide a characterization of the radii minimal projections of
  polytopes onto $j$-dimensional subspaces in 
  Euclidean space $\E^n$.  
  Applied on simplices this characterization allows to reduce
  the computation of an outer radius 
  to a computation in the circumscribing case or to the computation
  of an outer radius of a lower-dimensional simplex.
  In the second part of the paper, we use this characterization to 
  determine the sequence of outer $(n-1)$-radii of
  regular simplices (which are the radii of smallest enclosing
  cylinders).
  This settles a question which arose from the incidence that
  a paper by Wei{\ss}bach (1983) on this determination was erroneous.
  In the proof, we first reduce the problem to a constrained optimization 
  problem of symmetric polynomials and then to an optimization 
  problem in a fixed number of variables with additional integer 
  constraints. 
\end{abstract}

\maketitle


\section{Introduction}

Let 
$\mathcal{L}_{j,n}$ be the set of all $j$-dimensional linear
subspaces (hereafter \emph{$j$-spaces}) in $n$-dimensional
Euclidean space $\E^n$.
The \emph{outer $j$-radius} $R_j(C)$ 
of a convex body $C \subset \E^n$ is the 
radius of the smallest enclosing \emph{$j$-ball} ($j$-dimensional ball) 
in an optimal orthogonal projection of $C$ onto a $j$-space
$J \in \mathcal{L}_{j,n}$, where the optimization is
performed over $\mathcal{L}_{j,n}$.
The optimal projections are called \emph{$R_j$-minimal projections}. 
In this paper we show the following results:

\begin{thm} \label{th1}
 Let $1 \le j \le n < m$ and $P = \conv\{v^{(1)},\dots,v^{(m)}\}
 \subset \E^n$ be an $n$-polytope.
 Then one of the following is true.
 \begin{enumerate}[a)] 
  \item In every $R_j$-minimal projection
    of $P$ there exist $n+1$ affinely independent vertices of $P$ which
    are projected onto the minimal enclosing $j$-sphere.
  \item $j \ge 2$ and $R_j(P) = R_{j-1}(P \cap H)$ for some hyperplane
    $H = \aff \{v^{(i)}: i \in I\}$ with $I \subset \{1,\dots,m\}$.
 \end{enumerate}
If $j=1$ or if $P$ is a regular simplex then always case a) holds.
Moreover, the number $\nu$ of affinely independent vertices projected
onto the minimal enclosing $j$-sphere is at least $n-j+2$ and there
exists a ($\nu-1$)-flat $F$ such that $R_j(P) = R_{j+\nu-n-1}(P \cap
F)$.

\noindent
The bound $n-j+2$ is best possible.
\end{thm}

Theorem \ref{th1} allows to reduce the computation of an outer radius 
of a simplex to the computation in the circumscribing case
or to the computation of an outer radius of a 
facet of the simplex. 

Using this theorem, the second part of the paper shows
the following result on the outer $(n-1)$-radius, which is the 
radius of a smallest enclosing cylinder.

\begin{thm} \label{th2}
Let $n \ge 2$ and $T_1^n$ be a regular simplex in $\E^n$ with edge 
length 1. Then
\[
  R_{n-1}(T_1^n) \ = \left\{ \begin{array}{cl}
  \sqrt{\frac{n-1}{2(n+1)}} & \text{if $n$ is odd,} \\
  \frac{2n-1}{2\sqrt{2n(n+1)}} & \text{if $n$ is even.}
  \end{array} \right.
\]
\end{thm}

The case $n$ odd has already been settled independently by
Pukhov~\cite{pukhov-1980} and
Wei{\ss}bach~\cite{weissbach-83-correct}, who both left the even
case open in their papers. Pukhov's results also determine $R_j(T_1^n)$
for $j<n$. 
There also exists a later paper on $R_{n-1}(T^n)$ for even $n$
\cite{weissbach-83-erroneous}, but as pointed out
in \cite{brandenberg-theobald-2002} the proof contained a crucial
error.\footnote{After \cite{brandenberg-theobald-2002} had been
completed, 
Bernulf Wei{\ss}bach suggested to work jointly on a new
proof for the even case. Unfortunately, he died on 9th
June 2003.}
Thus Theorem~\ref{th2} (re-)completes the determination of the
sequence of outer
$j$-radii of regular simplices (see also \cite{brandenberg-regular-2002}).

Studying radii of polytopes is a fundamental topic in convex geometry.
Motivated by applications in
computer vision, robotics, computational biology, functional analysis, 
and statistics (see \cite{gritzmann-klee-mathprog-93} and the references 
therein) there has been much interest from the computational point of
view. See \cite{brandenberg-theobald-2002,dmpt-2002,ssty-2000} for exact
algebraic algorithms,
\cite{har-peled-varadarajan-04,vvz-2002,ye-zhang-2003} 
for approximation algorithms,
and \cite{bgkkls-2003,gritzmann-klee-mathprog-93} for the computational
complexity. 
Reductions of smallest enclosing cylinders to circumscribing 
cylinders are used in exact algorithms
as well as for complexity proofs
(see, e.g., \cite[Theorem~1]{brandenberg-theobald-2002} and
\cite[Theorems~5.3--5.5]{gritzmann-klee-mathprog-93}),
and have previously been given only for $j \in \{1, n\}$ as well as for 
dimension~3.
Theorem~\ref{th1} generalizes and unifies these results.

Here, we use Theorem~\ref{th1} to reduce the computation of the 
outer $(n-1)$-radius of a regular simplex to the following 
optimization problem of symmetric polynomials 
in $n$ variables:\begin{equation}
\label{eq:mainsystem}
 \begin{array}{lclcr}
  & \min & \sum \limits_{i=1}^{n+1}  s_i^4 && \\
  \text{s.t.} & & \sum \limits_{i=1}^{n+1}  s_i^3 & = & 0 \, , \\
  && \sum \limits_{i=1}^{n+1}  s_i^2 & = & 1 \, , \\
  && \sum \limits_{i=1}^{n+1}  s_i & = & 0 \, .
 \end{array}
\end{equation}
Based on exploiting the
symmetries, we solve \eqref{eq:mainsystem} for any $n$
by reducing it to an optimization
problem in six variables with additional integer constraints.

The paper is structured as follows. In Section~\ref{se:prelim}, we
introduce the necessary notation. Section~\ref{se:circum} gives the
proof of Theorem~\ref{th1}. Section~\ref{se:reductionalgebraic} contains
the derivation of the optimization 
problem~\eqref{eq:mainsystem} and the proof of Theorem~\ref{th2}.
In Section~\ref{se:realnullstellensatz} we analyze the
different difficulty of
the even and the odd case of~\eqref{eq:mainsystem} from the
viewpoint of the Positivstellensatz.

\section{Preliminaries\label{se:prelim}}
\label{reduce} 

Throughout the paper we work in Euclidean space $\E^n$, 
i.e., $\mathbb{R}^n$ with the usual scalar product 
$x \cdot y = \sum_{i=1}^n x_i y_i$ and norm $||x|| = (x \cdot
x)^{1/2}$.
$\B^n$ and $\mathbb{S}^{n-1}$ denote the (closed) unit ball
and unit sphere, respectively. 
For a set $A \subset \E^n$, the linear hull of $A$ 
is denoted by $\lin(A)$, the affine hull by $\aff(A)$, and
the convex hull by $\text{conv}(A)$.

A set $C \subset \E^n$ is called a \emph{body} if it is
compact, convex and contains an interior point. 
Accordingly, we always assume that a
polytope $P \subset \E^n$ is full-dimensional (unless otherwise stated).
Let $1 \le j \le n$. A \emph{$j$-flat} $F$ (an affine subspace of
dimension $j$) is \emph{perpendicular} to a hyperplane $H$ with
normal vector $h$ if $h$ and $F$ are parallel. 
For $p,p' \in \E^n$ and subspaces $E \in \mathcal{L}_{j,n}$, $E' \in \mathcal{L}_{j',n}$,
a $j$-flat $F = p + E$ and a $j'$-flat $F' = p' + E'$ 
are \emph{parallel} if $E \cup E' = \lin(E \cup E')$.

A \emph{$j$-cylinder} is a set of the form $J + \rho \B^n$ with
an $(n-j)$-flat $J$ and $\rho > 0$. 
For a body $C \subset \E^n$, the outer $j$-radius
$R_j(C)$ of $C$ (as defined in the introduction)
is also the radius $\rho$ of a smallest enclosing 
$j$-cylinder of $C$. 
It follows from a standard compactness argument that this minimal
radius is attained (see, e.g., \cite{gritzmann-klee-dcg-92}).
  Let $1 \le j \le k < n$.
  If $C' \subset \E^n$ is a compact, convex set whose affine
  hull $F$ is a $k$-flat then $R_{j}(C')$ denotes the radius of a 
  smallest enclosing $j$-cylinder $\mathcal{C'}$ relative to
  $F$, i.e., $\mathcal{C'} = J' +
  R_{j}(C') (\B^n \cap F)$ with a $(k-j)$-flat $J' \subset F$. 

A \emph{simplex} $S := \conv \{v^{(1)},\dots,v^{(n+1)}\}$
(with $v^{(1)},\dots,v^{(n+1)} \in \E^n$ affinely independent)
is \emph{regular} if all its 
vertices are equidistant. Whenever a statement is invariant under
orthogonal transformations and translations we denote by 
$T^n$ \emph{the} regular simplex in $\E^n$ with edge length $\sqrt{2}$.
The reason for the choice of $\sqrt{2}$ 
stems from the following embedding of $T^n$
into $\E^{n+1}$.
Let $\mathcal{H}^n_\alpha = \{x \in \E^{n+1} \, : \, \sum_{i=1}^{n+1} x_i =
\alpha\}$. Then the \emph{standard embedding} 
$\mathbf{T}^n$ of $T^n$ is defined by
\[
  \mathbf{T}^n \ := \ 
   \text{conv}\left\{ e^{(i)} \in \E^{n+1} \, : \, 1 \le i \le {n+1} \right\} 
  \subset
  \mathcal{H}^n_1 \, ,
\]
where $e^{(i)}$ denotes the $i$-th unit vector in $\E^{n+1}$.
By 
$\mathcal{S}^{n-1} := \mathbb{S}^{n} \cap \mathcal{H}^n_0$
we denote the set of unit vectors parallel to $\mathcal{H}_1^n$.

A $j$-cylinder $\mathcal{C}$
containing some simplex $S$ is called a 
\emph{circumscribing} $j$-cylinder of $S$ if all the
vertices of $S$ are contained in the boundary of $\mathcal{C}$.

\section{Minimal and circumscribing $j$-cylinders\label{se:circum}}

It is well known that the (unique) 
minimal enclosing ball $B$ (i.e., the minimal enclosing
$n$-cylinder) of a polytope $P \subset \E^n$ may
contain only few vertices of $P$ on its boundary (e.g., two diametral
vertices) \cite[p.~54]{bonnesen-fenchel-b34}. However,  
in cases where less than $n+1$ vertices of $P$ are contained in
the boundary of
$B$, it is easy to see that there exists a hyperplane $H$ such that $P
\cap \mathrm{bd}(B) \subset H$ and the center of $B$ is contained in
$H$. Then the smallest enclosing ball of $P$ and the smallest
enclosing ball of $P \cap H$ relative to $H$ have the same radius. 

In \cite[Theorem 1.9]{gritzmann-klee-dcg-92} the following
characterization for the minimal enclosing 1-cylinder (two parallel
hyperplanes defining the width of the polytope) is given:

\begin{prop} \label{prop-j1}
 Any minimal enclosing 1-cylinder of a  
 polytope $P \subset \E^n$ contains at least $n+1$ affinely
 independent vertices of $P$ on its boundary.
\end{prop}

We give a characterization of the possible
configurations of minimal enclosing 
$j$-cylinders of polytopes for arbitrary $j$, unifying and
generalizing the above statements.

\begin{lemma} \label{non-parallel-case}
 Let $P = \conv\{v^{(1)},\dots,v^{(m)}\}$ be a polytope in $\E^n$, $1 \le j \le n-1$,
 and $J$ be an $(n-j)$-flat such that $\mathcal{C} = J + R_j(P) \B^n$
 is a minimal enclosing $j$-cylinder of $P$. Then for every
 $I \subset \{1, \ldots, m\}$ such that
 $\{ i \, : \, v^{(i)} \in \mathrm{bd}(\mathcal{C}) \} \subset I$
 and $H_I := \mathrm{aff} \{v^{(i)} \, : \, i \in I \}$
 is of affine dimension $n-1$, $J$ is parallel to $H_I$.
\end{lemma}

\begin{proof}
 Suppose that there exists a hyperplane $H:=H_I$ of this type with $J$
 not parallel to $H$. Let $\bar{n}:= |\{v^{(i)} \in H \, : \, 1
 \le i \le m \}|$. Without loss of generality we can assume
  $H = \{x \in \E^n \, : \, x_n = 0\}$ and $I = \{v^{(1)}, \ldots, v^{(\bar{n})}\}$.
  Hence, $v^{(\bar{n}+1)}, \ldots, v^{(m)} \not\in H \cup \text{bd}(C)$. 
 
 First consider the case that $J$ is not perpendicular to
 $H$. Let $p,s^{(1)},\dots, $ $s^{(n-j)} \in \E^n$ such that $J = p +
 \lin \{ s^{(1)}, \ldots, s^{(n-j)} \}$. Since, by assumption, $J$ is
 not parallel to $H$, we can assume $p = 0 \in J \cap H$,
 $s_n^{(1)} = \dots = s_n^{(n-j-1)} = 0$ and $s_n^{(n-j)} > 0$.
 For every $s_n' \in (0, s_n^{(n-j)})$ and
 $s' := (s_1^{(n-j)}, \dots, s_{n-1}^{(n-j)}, s_n') \in \E^n$ let 
 $J' = p + \lin \{ s^{(1)}, \ldots, s^{(n-j-1)}, s' \}$. Geometrically,
 $J'$ results from $J$ by rotating $J$ towards the
 hyperplane $H$ in such a way that the orthogonal projection of $J$
 onto $H$ remains invariant (see Figure~\ref{fi:rotation1}). 

\ifpictures
\begin{figure}[h]
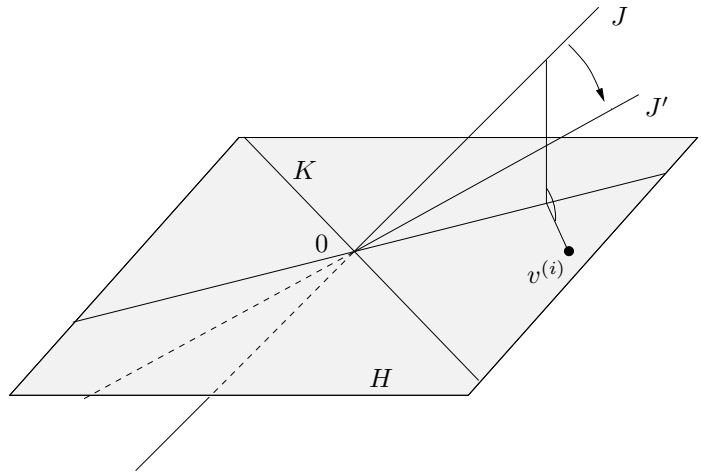

  \centerline{\mbox{} \input pictures/rotate2-new.pstex_t \mbox{}}
  \caption{\label{fi:rotation1} 
     For $n=3$ and $j=2$ the figure shows how the underlying flat $J$ 
     of the $j$-cylinder $\mathcal{C}$ is rotated towards its 
     orthogonal projection onto the plane $H$. The 
     distances between the vertices 
     $v^{(i)}$, $1 \le i \le \bar{n}$, and the $j$-cylinder axis are not
     increased, and decreased if $v^{(i)} \not\in K$.} 
  \end{figure}
\fi

 Since $J$ and $H$ are not perpendicular we obtain $J \neq J'$, and
 because $v^{(1)},\dots,v^{(\bar{n})} \in H$ that  
 \begin{equation}
  \label{dist}
  \mathrm{dist}(v^{(i)}, J') \le \mathrm{dist}(v^{(i)}, J) \, , \quad
  1 \le i \le \bar{n} \, ,
 \end{equation}
 where $\text{dist}(\cdot,\cdot)$ denotes the Euclidean distance.
 In~\eqref{dist}, ``$<$'' holds whenever $v^{(i)} \not\in K :=
 J^{\perp} \cap H$. Obviously, $\dim(K)=j-1$.
 If none of the $v^{(i)}$ lies in $K \cap \mathrm{bd}(\mathcal{C})$
 then, by choosing $s'_n$ sufficiently  
 close to $s_n^{(n-j)}$, all vertices of $P$ lie in the interior of
 $\mathcal{C}' = J' + R_j(P) \B^n$, a contradiction to the minimality
 of $\mathcal{C}$. Hence, there must be some vertex of $P$ in $K \cap
 \mathrm{bd}(\mathcal{C})$.  

 Let $\bar{k} := |\{ v^{(i)} \in K \cap \mathrm{bd}(\mathcal{C})  \, : \, 1 \le i \le m \} |$.
 By renumbering the vertices we can assume that
 $v^{(1)},\dots,v^{(\bar{k})} 
 \in K \cap \mathrm{bd}(\mathcal{C})$. Let $F := \conv
 \{v^{(1)},\dots,v^{(\bar{k})}\}$ and $k:= \dim F$. 

 Now assume $F \cap J = \emptyset$. We have shown above that for
 sufficiently small $s'_n$ the rotation from $J$ to $J'$ keeps all
 vertices within the $j$-cylinder $\mathcal{C}'$ and
 $v^{(1)},\dots,v^{(\bar{k})}$ are the only 
 vertices on $\mathrm{bd}(\mathcal{C}')$. Let $J''$ be a translate of
 $J'$ with $\text{dist}(J'',F) < \text{dist}(J',F)$, and $J''$
 sufficiently close to $J'$ to keep $v^{(\bar{k}+1)},\dots,v^{(m)}$ within
 the interior of $\mathcal{C}'' = J'' + R_j(P) \B^n$. Then all
 vertices of $P$ lie in the interior of $\mathcal{C}''$, again a contradiction.

 It follows that $F \cap J \neq \emptyset$, and since $F \subset K =
 J^\bot \cap H$ that $F \cap J = \{ p \} = \{ 0 \}$. 
 Since $\dist(p,v^{(i)}) = R_j(P)$ for all $i \in \{1,\dots,\bar{k}\}$ and
 since $p \in F$, it follows that $p$ is the unique 
 center of the smallest enclosing $k$-ball of $F$. 
 Now let $J'''$
 result from $J'$ by rotating $J'$ around the origin towards a
 direction in
 $\R^n \setminus (\bigcup_{i=1}^{\bar{k}} (v^{(i)})^{\perp})$.
 This set of directions is nonempty, since every subspace in the
 union is only of dimension $n-1$.
 For $i \in \{1,\dots,\bar{k}\}$ the property $\dist(v^{(i)},J) =
 \dist(v^{(i)},J') = \dist(v^{(i)},p)$ implies $\dist(v^{(i)},J''') <
 \dist(v^{(i)},J')$.
 Again, by keeping the rotation sufficiently small,
 $v^{(\bar{k}+1)},\dots,v^{(m)}$ remain in the interior of  
 $\mathcal{C}''' = J''' + R_j(P) \B^n$. Now, all vertices lie in the
 interior of $\mathcal{C}'''$, once more a contradiction. 
 
 Finally, consider the case that $J$ is perpendicular to $H$. Then $J \cap H$ 
 is an optimal $(n-1-j)$-flat for the $j$-radius of $P \cap H$
 (taken in $(n-1)$-dimensional space). However, it is easy to
 see that in this case any small perturbation $J'$ of $J$ with $J'
 \cap H = J \cap H$
 keeps $v^{(\bar{n}+1)},\dots,v^{(m)}$ within the $j$-cylinder, not
 increasing the distances of all the other vertices to the new
 $(n-j)$-flat. 
 Indeed, the case $J$ perpendicular to $H$ describes a local maximum.
 So the same argument as in the non-perpendicular
 case applies to show a contradiction.
\end{proof}

\begin{cor} \label{simplex-facet}
  Let $S=\conv\{v^{(1)},\dots,v^{(n+1)}\}$ be a simplex in $\E^n$ and
  $S^{(i)}$ be the facet of $S$ with $v^{(i)} \not\in
  S^{(i)}$. Let $1 \le j \le n-1$ and $J$ be an $(n-j)$-flat such that
  $\mathcal{C} = J + R_j(S) \B^n$ is a minimal enclosing $j$-cylinder
  of $S$. 
  Then $J$ is parallel to every $S^{(i)}$ for which $v^{(i)}
  \not\in \mathrm{bd}(\mathcal{C})$. 
\end{cor}

\begin{lemma} \label{parallel-case}
 Let $P = \conv\{v^{(1)},\dots,v^{(m)}\}$ be a polytope in $\E^n$, $1 \le j \le n$, and
 $J$ be an $(n-j)$-flat such that $\mathcal{C} = J + R_j(P) \B^n$ is a
 minimal enclosing $j$-cylinder of $P$. If there exists a hyperplane
 $H_I=\aff\{v^{(i)}:i \in I\}$ which is parallel to $J$, then one of
 the following holds:
 \begin{enumerate}[a)]
  \item There exists a vertex $v^{(i)} \not \in H_I$ that lies on the
    boundary of $\mathcal{C}$; or
  \item $j \ge 2$, $J \subset H_I$, and $R_j(P) = R_{j-1}(P \cap H_I)$.   
 \end{enumerate}
\end{lemma}

\begin{proof} By Proposition~\ref{prop-j1}, for $j=1$ always a) holds;
 so let $j \ge 2$, and suppose neither a) nor b) holds. 
 Since b) does not hold there exist $(n-j)$-flats parallel to $J$ and
 closer to $H_I$, and since a)
 does not hold, for any such $(n-j)$-flat $J'$, such that all vertices $v^{(i)} \not
 \in H_I$ stay within $\mathcal{C}$, the distances from the
 vertices $v^{(i)}$, $i \in 
 I$, to $J'$ are strictly smaller than their distances to $J$. Hence
 $\mathcal{C}$ cannot be a minimal enclosing cylinder.
\end{proof}

 In the case that $P$ is a simplex, the proof can be carried out more
 explicitly:

 Let $P^{(n+1)}$ be the facet of $P$ not including the vertex
 $v^{(n+1)}$. Suppose that $J$ is parallel to $P^{(n+1)}$, that
 $P^{(n+1)} \subset H := \{x \in \E^n : x_n = 0 \}$, and that
 $v_n^{(n+1)}>0$.  
 Let $p \in J$. Since $v_n^{(n+1)}>0$ it follows $p_n \ge 0$ and
 obviously 
 \begin{equation} \label{pn1} 
  R_j(P) \ge v_n^{(n+1)} - p_n .
 \end{equation}  
 On the other hand, since $J$ is parallel to $P^{(n+1)}$,
 \begin{equation} \label{pn2} 
  R_j(P)^2 = R_{j-1}^2(P^{(n+1)}) + p_n^2.
 \end{equation} 
 Let 
 \begin{equation} \label{p*}
  p_n^* = \frac{(v_n^{(n+1)})^2 - R_{j-1}^2(P^{(n+1)})}{2
 v_n^{(n+1)}}
 \end{equation}
 be the unique minimal solution for $p_n$ to
 \eqref{pn1} and \eqref{pn2}. Due to $p_n \ge 0$, we
 obtain $p_n = \max\{0,p_n^*\}$. Now, we see that case a) holds if
 $p_n = p_n^*$ and case b) if $p_n=0$. 

\bigskip

 Statements \ref{prop-j1}--\ref{parallel-case} almost complete the
 proof of Theorem \ref{th1}. 
 If the number $\nu$ of affinely independent vertices of $P$ lying on
 the boundary of $\mathcal{C}$ is at most $n$,
 it follows from
 Lemma \ref{non-parallel-case} and \ref{parallel-case} that case b) of
 Theorem \ref{th1} must hold. Moreover, if $\nu \le n-1$ we can again apply
 these lemmas on the 
 lower-dimensional polytope $P \cap H_I$ with $H_I$ as
 in Lemma \ref{parallel-case}. 
 Now we can iterate this argument. If during this iteration
 the outer 1-radius of a polytope $P'$ has to be
 computed, then by Proposition~\ref{prop-j1} the
 minimal enclosing 1-cylinder touches at least
 $\dim(P')+1$ affinely independent vertices. From the same iterative
 argument it follows that $R_j(P) = R_{j+\nu-n-1}(P \cap F)$ for some
 ($\nu-1$)-flat $F$.

 Suppose $S=\conv\{v^{(1)},\dots,v^{(n+1)}\}$ is a simplex in $\E^n$,
 and $\bar{J}$ an $(n-j)$-flat, such that
 \begin{eqnarray*}
  \dist(v^{(1)},J)&=& \cdots \ = \ \dist(v^{(n-j+2)},J) \ = \
  R_1(\conv\{v^{(1)},\dots,v^{(n-j+2)}\})  \\
  &>& \dist(v^{(n-j+3)},J) \ \ge \ \cdots \ \ge \ \dist(v^{(n+1)},J).
 \end{eqnarray*}
 Then obviously $R_j(S)=R_1(\conv\{v^{(1)},\dots,v^{(n-j+2)}\})$ and
   at most $n-j+2$ vertices are situated on the boundary of the
   minimal enclosing $j$-cylinder.

\bigskip

 The last point which remains to show is
 that every minimal enclosing $j$-cylinder of the regular simplex $T^n$
 also circumscribes $T^n$. Because of Proposition \ref{simplex-facet}
 it suffices to show that the value $p_n^*$ in~\eqref{p*} is positive
 for all $1 \le j \le n-1$, showing that b) in Lemma~\ref{parallel-case}
 never holds for $T^n$. 
 
 However, in almost all cases the desired circumscribing property follows
 already from \cite{pukhov-1980}, see also \cite{brandenberg-regular-2002}. 

\begin{prop} \label{known}
 For $1 \le j \le n$ it holds $R_j(T^n) \ge \sqrt{\frac{j}{n+1}}$.
 If $n$ is odd or $j \not \in \{1,n-1\}$, then $R_{j}(T^n) =
 \sqrt{\frac{j}{n+1}}$, and every minimal enclosing $j$-cylinder of
 $T^n$ is a circumscribing  
 $j$-cylinder of $T^n$.
\end{prop}

 We can easily apply Proposition \ref{known} to compute $p_n^*$
 if $n$ is even.

 \begin{lemma} \label{even-case}
  Let $1 \le j \le n-1$. If $S=T^n$ then always case a) in Theorem
  $\ref{th1}$ holds.
 \end{lemma}

 \begin{proof} We can assume $j \ge 2$, since otherwise b) cannot hold.
  We use the notation as in Lemma \ref{parallel-case}.  Because of
  Proposition \ref{known} it suffices to consider the case where $n$
  is even, and as mentioned 
  above the proof is complete if we show that $p_n^*$ is positive. 
  Since Proposition~\ref{known} yields $R_{n-1}(T^{n-1}) =
  \sqrt{(n-1)/n}$, we have  $v^{(n+1)}_n = \sqrt{2-(n-1)/n} =
  \sqrt{(n+1)/n}$. Also by
  Proposition~\ref{known}, $R_{j-1}(T^{n-1}) = \sqrt{(j-1)/n}$ and therefore 
  $$p_n^* \ = \ \frac{n-j+2}{2\sqrt{n(n+1)}} \ > \ 0 \, .$$
 \end{proof}
  
 Choosing an optimal $(n-1)$-cylinder among those 
 parallel to a facet of $T^n$ with
 $p_n^* = \frac{3}{2\sqrt{n(n+1)}}$, gives an upper 
 bound for the outer $(n-1)$-radius of a regular simplex,
 $$ R_{n-1}(T^n) \  \le  \ \frac{2n-1}{2\sqrt{n(n+1)}} \, .$$
 Theorem \ref{th2} states that for even $n$ this bound is tight.

\section{Reduction to an algebraic optimization problem\label{se:reductionalgebraic}}

In this section, we provide an algebraic formulation for a minimal
circumscribing $j$-cylinder $J + \rho (\B^{n+1} \cap \mathcal{H}_0^n)$ of
the regular simplex $\mathbf{T}^n$ 
in standard embedding. Let $J = p + \lin\{s^{(1)}, \ldots,
s^{(n-j)}\}$ with pairwise orthogonal $s^{(1)},\dots,s^{(n-j)} \in
\mathcal{S}^{n-1}$, and $p$ be contained in the orthogonal 
complement of $\lin \{s^{(1)},\dots,s^{(n-j)}\}$.
The orthogonal projection $P$
of a vector $z \in \mathcal{H}^{n}_{1}$ onto the orthogonal
complement of $\lin \{s^{(1)},\dots,s^{(n-j)}\}$
(relative to $\mathcal{H}^{n}_{1}$) 
can be written as  
$P( z) = (I - \sum_{k=1}^{n-j} s^{(k)} (s^{(k)})^T) z$, where
$I$ denotes the identity matrix. 
Hence, for a general polytope with vertices $v^{(1)},\dots,v^{(m)}$ 
(embedded in $\mathcal{H}^{n}_{1}$) the
computation of the square of $R_j$ 
can be expressed by the following optimization
problem.
Here, we use the convention $x^2:=x \cdot x$.
 
\begin{equation} \label{gen-poly}
  \begin{array}{l@{\qquad \qquad}rrclcr}
   && \min \, \rho^2 \\
   \text{(i)} & \text{s.t.} 
    & (p - P v^{(i)})^2 & \le & \rho^2
    \, , & \quad i=1,\dots,m, \\ 
    \text{(ii)} & & p \cdot s^{(k)} & = & 0 \, , & k=1,\dots,n-j,\\
   \text{(iii)} & & s^{(1)}, \ldots, s^{(n-j)} & \in & \mathcal{S}^{n-1}, & \text{pairwise
   orthogonal,}\\ 
   \text{(iv)} & & p & \in & \mathcal{H}_{1}^{n} \, .
  \end{array}
 \end{equation}
In the case of $\mathbf{T}^n$, (i) can be replaced by 
\begin{equation}
  \label{rhosquare}
  \text{(i')} \qquad \qquad
   \left(p - e^{(i)} + \sum \limits_{k=1}^{n-j} s_i^{(k)}
   s^{(k)}\right)^2 \ = \ \rho^2 \, ,  \quad i=1,\dots,n+1 \, , 
\end{equation}
where the equality sign follows from Theorem \ref{th1}. 
By (ii) and $s^{(k)} \in \mathcal{S}^{n-1}$, 
(i') can be simplified to 
\[
  \text{(i'')} \qquad \qquad
        p^2 - \rho^2 \ = \ \sum \limits_{k=1}^{n-j} (s_i^{(k)})^2 +
        2p_i - 1 \, ,  \quad i = 1, \ldots, n+1 \, .
\]
Summing over all $i$ gives $(n+1) (p^2 - \rho^2) = (n-j)+ 2-(n+1)$, i.e., 
$p^2 - \rho^2 = \frac{1-j}{n+1}$. We substitute this value 
into (i'') and obtain
$ p_i = \frac{1}{2}\left(\frac{n-j+2}{n+1} - \sum_{k=1}^{n-j}
(s_i^{(k)})^2\right)$. Hence,  
all the $p_i$ can be replaced in terms of the $s_i^{(k)}$,

\begin{eqnarray}
  \rho^2  & = & \frac{(2+(n-j))(2-(n-j))}{4(n+1)} + \frac{1}{4}
  \sum_{i=1}^{n+1} \left( \sum_{k=1}^{n-j} (s_i^{(k)})^2 \right)^2 +
  \frac{j-1}{n+1} \, , \label{transform} \\  
  p \cdot s^{(k)} & = & - \frac{1}{2} \sum_{i=1}^{n+1}
  \sum_{k'=1}^{n-j} (s_i^{(k')})^2 s_i^{(k)} \, . \nonumber
\end{eqnarray}

We arrive at the following characterization of the 
minimal enclosing $j$-cylinders:

\begin{thm} \label{th:polyformulation}
Let $1 \le j \le n$.
A set of vectors $s^{(1)},\dots,s^{(n-j)} \in \mathcal{S}^{n-1}$ spans
the underlying $(n-j)$-dimensional subspace of a minimal enclosing 
$j$-cylinder
of $\mathbf{T}^n \subset \mathcal{H}^n_1$ if and only if it is an 
optimal solution of the problem
\begin{equation}\label{eq:regsimplex0}
  \begin{array}{rrrcll}
   & \min & \sum \limits_{i=1}^{n+1} \left( \sum\limits_{k=1}^{n-j}
   (s_i^{(k)})^2 \right)^2 &\\ 
   \mathrm{s.t.} & 
      & \sum \limits_{i=1}^{n+1} \sum \limits_{k'=1}^{n-j} (s_i^{(k')})^2
   s_i^{(k)} & = & 0, & k=1,\dots,n-j  \, ,\\ 
      && s^{(1)}, \ldots, s^{(n-j)} & \in & \mathcal{S}^{n-1} &
   \text{pairwise orthogonal.} 
  \end{array}
 \end{equation}
\end{thm}

It is easy to see that in case $j=n-1$ 
the program~\eqref{eq:regsimplex0} reduces to~\eqref{eq:mainsystem}
stated in the introduction. 

By~\eqref{transform}, in order to prove
$R_{n-1}(T^n) = (2n-1)/(2\sqrt{n(n+1)})$ for even $n$, we have
to show that the optimal value of~\eqref{eq:mainsystem} is $1/n$.
We apply the following statement from 
\cite{brandenberg-theobald-2002}.

\begin{prop} \label{only3}
 Let $n \ge 2$. The direction vector $(s_1, \dots, s_{n+1})^T$ of any
 extreme circumscribing 
$(n-1)$-cylinder of $\mathbf{T}^n$ satisfies 
 $|\{ s_1, \ldots, s_{n+1} \}| \le 3$.
\end{prop}

For completeness we repeat the short proof.

\begin{proof}
We can assume $n \ge 3$. Let $s \in \mathcal{S}^{n-1}$ be the axis
direction of a locally extreme circumscribing $(n-1)$-cylinder.
Let $f(s) := \sum_{i=1}^{n+1} s_i^4$ be the objective function 
from~\eqref{eq:mainsystem}, let $g_1(s) := \sum_{i=1}^{n+1} s_i^3$,
$g_2(s) := \sum_{i=1}^{n+1} s_i^2 - 1$, and
$g_3(s) := \sum_{i=1}^{n+1} s_i$.
A necessary condition for a local extremum is that for any pairwise different
indices $a,b,c,d \in \{1, \ldots, n+1\}$,
\[
  \det \left( \begin{array}{cccc}
    -\frac{\partial f}{\partial s_a} &
    \frac{\partial g_1}{\partial s_a} &
    \frac{\partial g_2}{\partial s_a} &
    \frac{\partial g_3}{\partial s_a} \\ [0.4ex]
    -\frac{\partial f}{\partial s_b} &
    \frac{\partial g_1}{\partial s_b} &
    \frac{\partial g_2}{\partial s_b} &
    \frac{\partial g_3}{\partial s_b} \\ [0.4ex]
    -\frac{\partial f}{\partial s_c} &
    \frac{\partial g_1}{\partial s_c} &
    \frac{\partial g_2}{\partial s_c} &
    \frac{\partial g_3}{\partial s_c} \\ [0.4ex]
    -\frac{\partial f}{\partial s_d} &
    \frac{\partial g_1}{\partial s_d} &
    \frac{\partial g_2}{\partial s_d} &
    \frac{\partial g_3}{\partial s_d} \\
  \end{array} \right) \ = \
  -24 \, \det \left( \begin{array}{cccc}
    s_a^3 & s_a^2 & s_a & 1 \\
    s_b^3 & s_b^2 & s_b & 1 \\
    s_c^3 & s_c^2 & s_c & 1 \\
    s_d^3 & s_d^2 & s_d & 1 \\
  \end{array} \right) \ = \ 0 \, .
\]
The latter is a Vandermonde determinant,
which implies $|\{s_a,s_b,s_c,s_d\}| \le 3$.
\end{proof}

Using Proposition~\ref{only3}, \eqref{eq:mainsystem}
can be written as the following polynomial optimization problem
in six variables
with additional \emph{integer} conditions.
\begin{equation}\label{eq:regsimplex1}
  \begin{array}{l@{\qquad}rrcl}
   && \min \, k_1 s_1^4 + k_2 s_2^4 + k_3 s_3^4 \\
      \text{(i)} & \text{s.t.} & k_1 s_1^3 + k_2 s_2^3 + k_3 s_3^3 & = & 0 \, , \\
      \text{(ii)} && k_1 s_1^2 + k_2 s_2^2 + k_3 s_3^2 & = & 1 \, , \\
      \text{(iii)} && k_1 s_1 + k_2 s_2 + k_3 s_3 & = & 0 \, , \\
      \text{(iv)} && k_1 + k_2 + k_3 & = & n+1 \, , \\
      & \multicolumn{4}{r}{s_1, s_2, s_3 \in \mathbb{R}, \quad k_1, k_2, k_3 \in \mathbb{N}_0 \, .}
  \end{array}
 \end{equation}

Since the odd case of Theorem~\ref{th2} is 
well-known \cite{pukhov-1980,weissbach-83-correct},
we assume from now on that $n$ is even. The mindful reader
will notice that for odd $n$ the optimal value
of~\eqref{eq:regsimplex1} 
coincides with the optimal value of the real relaxation 
(where the condition $k_1, k_2, k_3 \in \mathbb{N}_0$ is replaced by
$k_1,k_2,k_3 \ge 0$). 

For $k_3 = 0$ the equality constraints in~\eqref{eq:regsimplex1} 
immediately yield $k_1 = k_2 = (n+1)/2 \not\in \N$, and
similarly, for $s_2 = s_3$ we obtain 
$k_1 = k_2 + k_3 = (n+1)/2 \not\in \N$. Hence, we can
assume that $s_1$, $s_2$, and $s_3$ are distinct and $k_1,k_2,k_3 \ge 1$.
Moreover, for $s_3 = 0$ the resulting optimal value is $1/n$
which will turn out to be the optimal solution.
Finally, by~(iii), not all of the $s_i$ have the same sign. Hence it
suffices to show that for $s_1 < 0$ and $s_3 > s_2 > 0$ every
admissible solution to the constraints  
of~\eqref{eq:regsimplex1} has value at least $1/n$.

The linear system of equations in $k_1,k_2,k_3$
defined by (i), (ii), and (iii) is regular, which is easily
seen from the Vandermonde computation
\[
  \det \left( \begin{array}{ccc}
    s_1^3 & s_2^3 & s_3^3 \\
    s_1^2 & s_2^2 & s_3^2 \\
    s_1   & s_2   & s_3 
  \end{array} \right) \ 
  = s_1 s_2 s_3 (s_1 - s_2)(s_1 - s_3)(s_2 - s_3) \neq 0 \, .
\]
Solving for $k_1,k_2,k_3$ yields
\begin{eqnarray}
  k_1 & = & \frac{s_2 + s_3}{-s_1(s_2 - s_1)(s_3 - s_1)} \, , \label{k1} \\
  k_2 & = & \frac{s_1 + s_3}{s_2(s_2 - s_1)(s_3-s_2)} \, , \label{k2} \\
  k_3 & = & \frac{-(s_1 + s_2)}{s_3(s_3 - s_1)(s_3-s_2)} \, . \label{k3}
\end{eqnarray}
Since all factors in the denominators are
strictly positive, \eqref{k2} and \eqref{k3} imply in particular 
$s_1 + s_3 > 0$ and $s_1 + s_2 < 0$.

Using (iv) in \eqref{eq:regsimplex1} we can express one of the $s_i$
by the others. Solving for $s_2$ gives
\begin{equation}
  \label{s2}
  s_2 \ = \  - \frac{s_1 + s_3}{(n+1) s_1 s_3 + 1} \, .
\end{equation}

Note that the denominator in $\eqref{s2}$ is
strictly negative. 

Our main goal now is to show $k_1 < (n+1)/2$
and then to use the integer condition to deduce $k_1 \le n/2$.
In order to achieve this, substitute \eqref{s2} into the inequality
$s_1+s_2 < 0$
which (in connection with $s_1+s_3 > 0$) allows
to conclude $s_3^2 > s_1^2 > 1/(n+1)$. 
Then substitute~\eqref{s2} into~\eqref{k1} which yields
\begin{align*}
  k_1 - \frac{n+1}{2} & =  
    - \frac{((n+1)s_1^2 - 1)((n+1)s_3(s_3-s_1) - 2)}
           {2 (s_3 - s_1)((n+1)s_1^2s_3 + 2s_1 + s_3)} \\
    & =  - \frac{((n+1)s_1^2 - 1)(((n+1)s_3^2-1) - ((n+1)s_1s_3+1))}
           {2 (s_3 - s_1)(((n+1)s_1 s_3 + 1) s_1 + (s_1 + s_3))} \ < \ 0 \, ,
\end{align*}
since all factors within the last fraction are positive.
Hence, $k_1 < (n+1)/2$ and since it is an integer $k_1 \le n/2$.
Similarly, although we do not need this, one can show 
$k_3 \le n/2$. However, note
that this bound does not hold for $k_2$.

Now it follows from $k_1 \le n/2$ and \eqref{k1} that
\begin{eqnarray*}
  0 & \le & 2 - n s_3^2 - 2 s_3^2 - n^2 s_1^3s_3 + n^2 s_1^2 s_3^2
            - n s_1^3 s_3 + n s_1^2 s_3^2 - 2 n s_1^2 + n s_1 s_3 \\
    & = & - 2 (n s_1^2 -1)((n+1)s_1 s_3 + 1) 
       + s_3 (s_1 + s_3) (n(n+1)s_1^2 - n -2) \, .
\end{eqnarray*}
That means at least one of the two terms of the sum must be non-negative
which gives $s_1 \le - 1/\sqrt{n}$. 
Moreover, $s_1 + s_3 > 0$ also implies $s_3 \ge 1/\sqrt{n}$.

Finally, we show that for any admissible solution to the
constraints of~\eqref{eq:regsimplex1}
the objective value is at least $1/n$. 
Replacing $k_1,k_2,k_3$ and $s_2$ in the objective function 
via \eqref{k1}--\eqref{s2} and using the inequalities $-s_1 < s_3$ and
$s_1 \le -1/\sqrt{n}$ obtained above we get   
\begin{eqnarray*}
  k_1 s_1^4 + k_2 s_2^4 + k_3 s_3^4 & = &  
  \frac{1}{n+1} + \frac{((n+1)s_1^2-1)((n+1)s_3^2-1)}{(n+1)(-(n+1)s_1s_3-1)} \\
  & \ge & \frac{1}{n+1} + \frac{((n+1) s_1^2 -1)}{n+1} \ \ge \ 
  \ \frac{1}{n+1} + \frac{1}{n(n+1)} \ = \ \frac{1}{n} \, .
\end{eqnarray*}
Hence, the optimal value 
of~\eqref{eq:regsimplex1} is $1/n$, and by our remark 
before Proposition~\ref{only3}
this completes the proof of
Theorem \ref{th2}.

\section{Connections to the Positivstellensatz\label{se:realnullstellensatz}}

We close the paper by discussing the greater difficulty of the even
case of computing $R_{n-1}(T^n)$
compared to the odd case, by analyzing problem~\eqref{eq:mainsystem} 
from the viewpoint of the Positivstellensatz \cite{stengle-74}.
This theorem states the existence of a certificate whenever a system
of polynomial equalities and inequalities does not have a solution, 
and it can be regarded as a common generalization of 
Hilbert's Nullstellensatz and of linear programming duality.
For our purposes, it suffices to consider the following 
version of Putinar (see \cite{prestel-delzell-2001, schweighofer-2004}).
For $n \in \N$ let $\R[x] = \R[x_1, \ldots, x_n]$ denote the ring of
polynomials in $x_1, \ldots, x_n$, and let 
\[
 \sum \R[x]^2 \ = \ \left\{ \sum_{j=1}^k b_j^2 \text{ for some } k \in \N,
  b_1, \ldots, b_k \in \R[x] \right\}
\]
be the set of all finite sums of squares of polynomials. Set $g_0 := 1$.

\begin{prop} If the polynomials $f,g_1, \ldots, g_m \in \R[x]$ 
satisfy
\[
  f(x) > 0 \quad \text{ for all } x \in S:= \{x \in \R^n \, : \ 
           g_1(x) \ge 0, \ldots, g_m(x) \ge 0 \}
\]
and
\begin{equation}
  \label{eq:quadmodule}
  M \ := \ \left\{ \sum_{i=0}^m \sigma_i g_i \, : \, 
   \sigma_0, \ldots, \sigma_m \in \sum \R[x]^2 \right\}
\end{equation}
contains the polynomial $1-\sum_{i=1}^n x_i^2$, then $f \in M$.
\end{prop}

Hence, in order to prove that the optimal value 
of~\eqref{eq:mainsystem} is bounded from below by
some given value $\alpha$,
it suffices (by compactness of the admissible set)
to show the existence of such a representation for 
$f(s) := \sum_{i=1}^{n+1} s_i^4 - \alpha + \varepsilon$
in terms of  
$g_1(s) := \sum_{i=1}^{n+1} s_i^3$, 
$g_2(s) := - \sum_{i=1}^{n+1} s_i^3$, 
$g_3(s) := \sum_{i=1}^{n+1} s_i^2 -1$,
$g_4(s) := - \sum_{i=1}^{n+1} s_i^2 +1$,
$g_5(s) := \sum_{i=1}^{n+1} s_i$,
$g_6(s) := - \sum_{i=1}^{n+1} s_i$
for every $\varepsilon > 0$.

Bounding the degrees of the polynomials $\sigma_i g_i$ by a fixed constant
in~\eqref{eq:quadmodule} serves to give lower bounds on the
minimum. These relaxations can be computed by semidefinite programming (SDP)
and are at the heart of current developments in SDP-based constrained
polynomial
optimization (see~\cite{lasserre-2001,parrilo-2003}).

For the case $n$ odd of~\eqref{eq:mainsystem} there exists
a simple polynomial identity
\begin{equation}
  \label{identity}
  \sum_{i=1}^{n+1} s_i^4 - \frac{1}{n+1} 
  \ = \ 
  \frac{2}{n+1} \left( \sum_{i=1}^{n+1} s_i^2 - 1 \right)
  + \sum_{i=1}^{n+1} \left( s_i^2 - \frac{1}{n+1} \right)^2
\end{equation}
which shows that the minimum is bounded from below by $1/(n+1)$,
and since this value can be attained
by $s_1 = \ldots = s_{(n+1)/2} = - s_{(n+3)/2} = \ldots = 
  - s_{n+1} = 1/{\sqrt{n+1}}$,
the minimum is $1/(n+1)$. For any $\varepsilon > 0$, adding $\varepsilon$
on both sides of~\eqref{identity} yields a representation of the 
positive polynomial on the left side as a sum of squares of the $g_i$.
Note that for every odd $n$ this representation 
uses only polynomials $\sigma_i g_i$ of (total) degree at most $4$. 

For the case $n$ even (with minimum $1/n$)
the situation is quite different. 
A computer calculation using the Software~{\sc GloptiPoly}
\cite{gloptipoly}
showed that already for $n=4$ it is necessary to go up to degree 8 
to find the Positivstellensatz-type certificate of optimality.
Since from a practical point of view the computational efforts 
of this calculation drastically increase with the number of variables, 
we do not know up to which degree it is necessary to go for $n=6$.

\subsection*{Acknowledgments.} We would like to thank Keith Ball for
pointing out a substantial simplification in the proof of Theorem \ref{th2}
and Pablo Parrilo for helpful discussions.

\providecommand{\bysame}{\leavevmode\hbox to3em{\hrulefill}\thinspace}
\providecommand{\MR}{\relax\ifhmode\unskip\space\fi MR }
\providecommand{\MRhref}[2]{%
  \href{http://www.ams.org/mathscinet-getitem?mr=#1}{#2}
}
\providecommand{\href}[2]{#2}

\bigskip

 
\end{document}